# On the Theory of Colorful Graphs

Dhananjay P. Mehendale
Sir Parashurambhau College, Tilak Road, Pune-411030,
India

## Abstract

The theory of colorful graphs can be developed by working in Galois field modulo ($p$), $p > 2$ and a prime number. The paper proposes a program of possible conversion of graph theory into a pleasant colorful appearance. We propose to paint the usual black (indicating presence of an edge) and white (indicating absence of an edge) edges of graphs using multitude of colors and study their properties. All colorful graphs considered here are simple, i.e. not having any multiple edges or self-loops. This paper is an invitation to the program of generalizing usual graph theory in this direction.

**1. Introduction:** In the usual graph theory one works in Galois field modulo (2), $GF(2)$. In the theory of **colorful graphs** (**cgraphs**, hereafter) we work in Galois field modulo ($p$), $GF(p)$, where $p > 2$ is some suitable prime (if one doesn't need to use the properties of field then any positive integer greater than two will do).

Galois field mod ($p$), $p$ prime, is a finite field containing elements {0, 1, 2, **…,** $p-1$} with modulo $p$ addition and multiplication. We paint the edges of a cgraph with colors from the color set $C$ = $\{C_0, C_1, \cdots, C_{p-1}\}$ and use the numbers (elements of the Galois field) {0, 1, 2, **…,** $p-1$} as matrix elements in their associated matrices, for showing the presence of these edges with a particular color in the cgraph under consideration. Note that by color $C_0$ (say, white) in a cgraph we represent the **absence** of an edge as is done in case of graphs.

**2. Multicolored Graphs:** If we paint the edges of usual simple graph and convert it into a cgraph using certain colors of our choice and decide to treat any colored edge as an edge in the usual graph theory and an absence of an edge (i.e. a white colored edge) as an absence in usual sense then we will have one kind of development in which all the definitions and propositions in usual graph theory will survive with only change where edges get replaced by colored edges: For example, we will get new straightforwardly modified definitions and results like

**Definition 2.1:** A cgraph is said to be **connected** if there exists a multicolored path joining any two vertices.

**Proposition 2.1:** A cgraph is disconnected if and only if its vertex set can be partitioned into (at least) two (or more) nonempty disjoint subsets such that there exists no colored edge (except white edge which color is taken to denote absence of an edge) joining any two vertices belonging to the disjoint vertex sets.

**Definition 2.2:** The **degree** of a vertex is the count of colored edges incident on it.

**Proposition 2.2:** If a cgraph (connected or disconnected) has exactly two vertices of odd degree, there must be a path of multicolored edges joining these two vertices.

**Proposition 2.3:** A simple cgraph with $n$ vertices and $k$ components can have at most $(n-k)(n-k+1)/2$ colored edges.

etc. etc.

Is it possible to have some interesting results that one gets by just replacing edges by colored edges? Here I mention three different application avenues of taking multicolored edges in usual graphs:

**2.4 Tight Packing of Complete Cgraphs:** Refer to theorem 6.1 in [1]. There, it is shown that the for the existence of a finite projective plane of order n and the existence of tight packing of certain complete graphs (distinguished by different colors) into a complete graph of certain large order are equivalent things. Let us consider the following special cases:

**Example 2.4.1:** Representation for Projective plane of order $n = 2$:

The cgraph given below, representing projective plane of order two, contains, in all, 35 triangles out of which only 7 triangles have same colors to all of their respective edges while the other 28 triangles have different colors to all of their respective edges!!

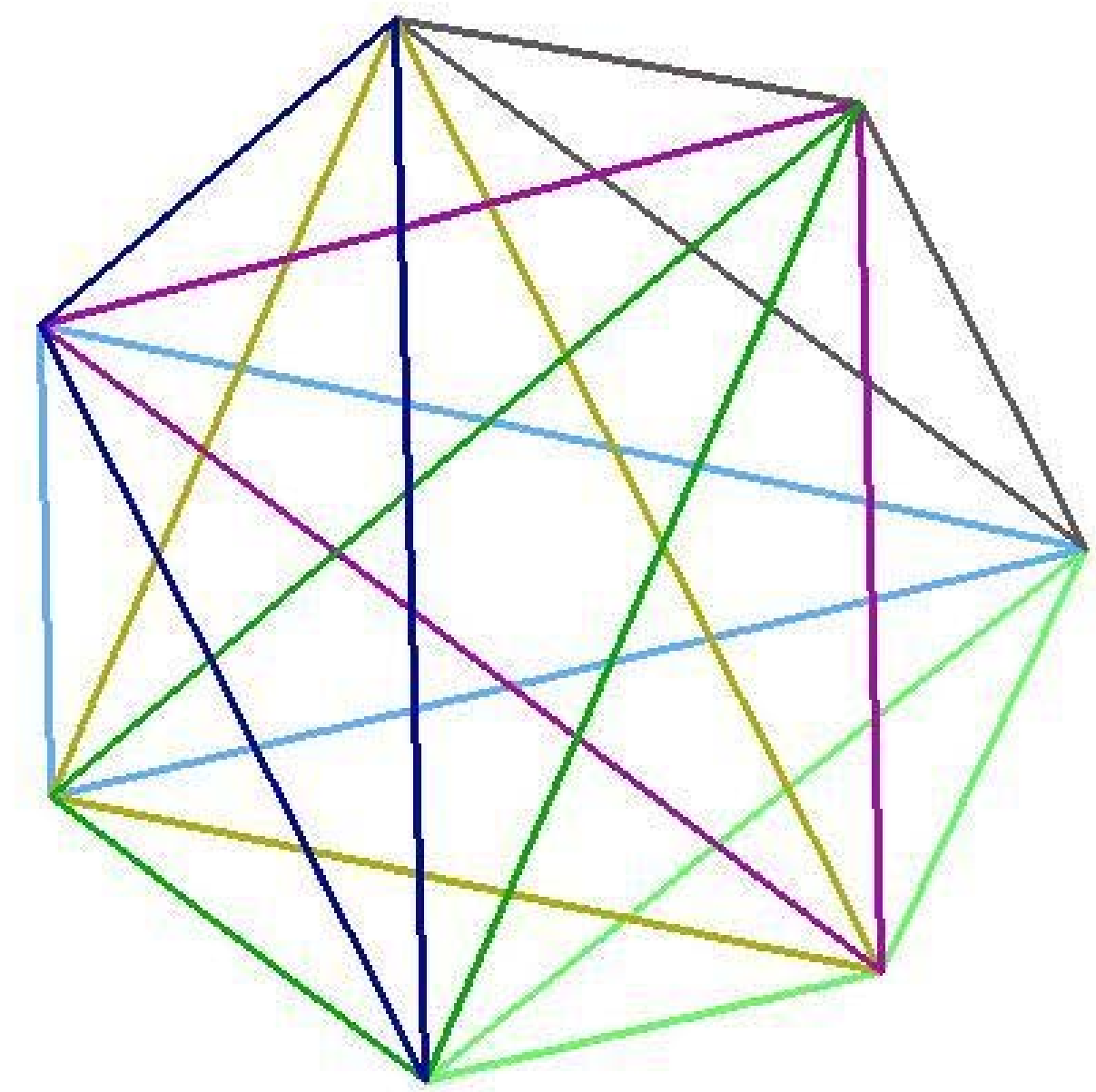

**Example 2.4.2:** Representation for Projective plane of order $n = 3$:

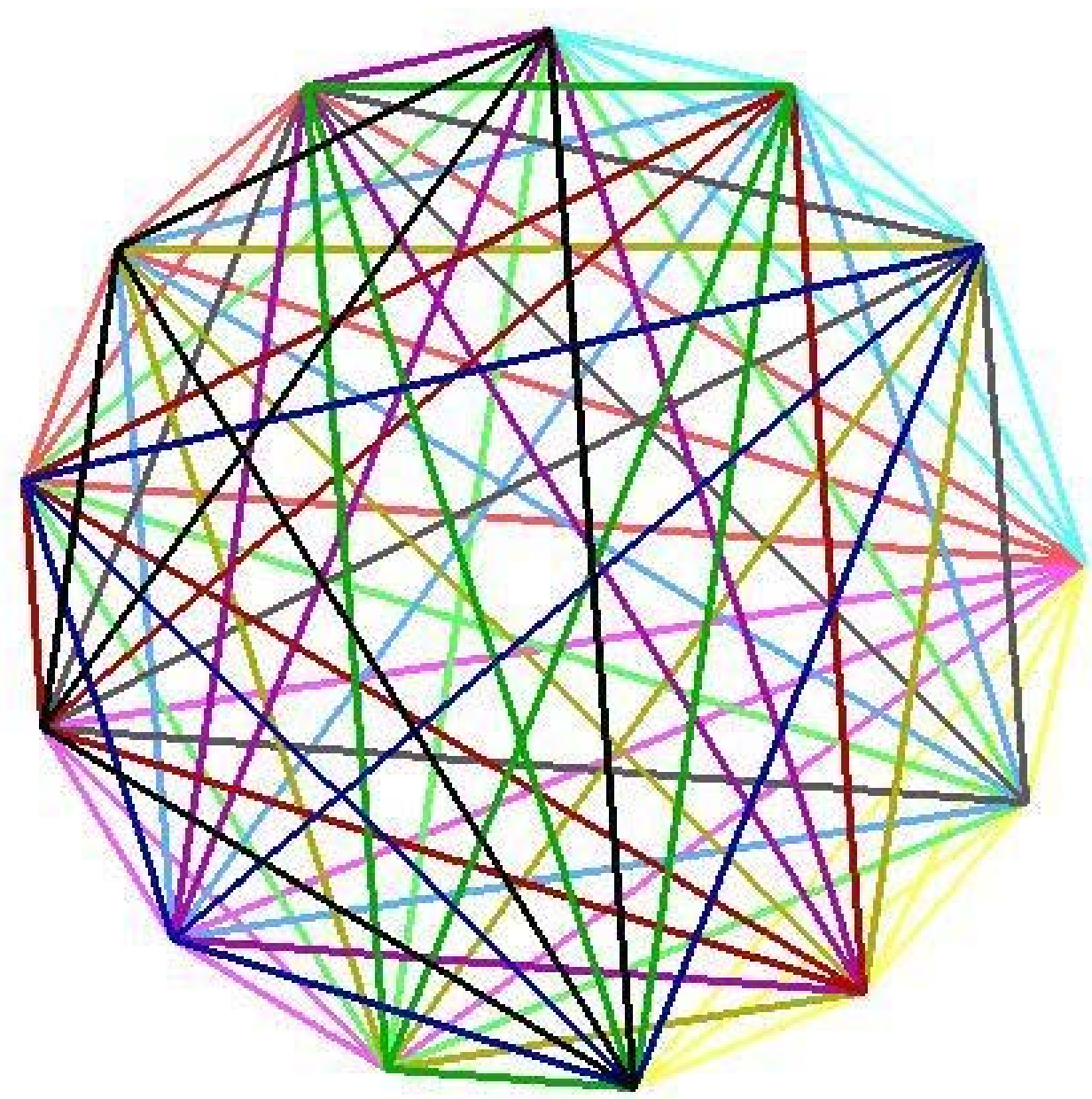

In the above example there are exactly 13 complete monochromatic cgraphs on 4 points, each with different color,  having same color to their respective edges!!

**2.5 An Application to Job Assignment Problem:** In order to illustrate the use of a (bipartite) cgraph we consider its simple application to job assignment problem:
Suppose we have *m* jobs and *n* people, and each person can do some of the jobs. Can we make assignments to fill the jobs? We model the available assignments by colored edges from vertex set $\{j_1, j_2, \cdots, j_m\}$ representing jobs to vertex set $\{p_1, p_2, \cdots, p_n\}$ representing persons in a cgraph. All edges emerging from vertex $j_i$, will be colored with color $i$, $1 \le i \le m$. We then construct a matrix M = $[a_{uv}]$ of size $(n \times m)$ such that $a_{uv} = v$, if there is an edge from vertex $p_u$ to vertex $j_v$, and $a_{uv}$ = 0, otherwise. Then, every nonzero determinantal monomial represents a proper job assignment.

**Example 2.5.1:** Consider the following matrix *M* (in a tabular form) corresponding to a cgraph:

$$M = \begin{array}{c|cccc} & j_1 & j_2 & j_3 & j_4 \\ \hline p_1 & 1 & 2 & 0 & 0 \\ p_2 & 1 & 0 & 3 & 0 \\ p_3 & 0 & 2 & 0 & 4 \\ p_4 & 0 & 0 & 0 & 4 \end{array}$$

The proper job assignment (from a nonzero determinantal monomial) is:

| Job | Person |
|---|---|
| $p_1$ | $j_1$ |
| $p_2$ | $j_3$ |
| $p_3$ | $j_2$ |
| $p_4$ | $j_4$ |

**2.6 Polya theory for cgraphs:** Polya's counting theorem is the most powerful tool in graph enumeration. Here we see that it can also work with the same effectiveness for cgraph enumeration problems. We illustrate here its use for counting (simple) cgraphs.

It is clear to see that the counting of cgraphs having some predefined property can be achieved by just replacing the figure counting series, which now change, and which are to be put in the place of the variables in the cycle index of the associated permutation group.

**2.6.1 Enumeration of Simple Cgraphs with Polya's Theorem:** A cgraph defined over the Galois field modulo$(p)$ can have $p$ colors for its edges from the color set $C = \{C_0, C_1, \cdots, C_{p-1}\}$, where color $C_0$ represents the white color (or the absence of a edge in the usual sense). The permutation group that is relevant in the case of graph counting remains the same for cgraph counting. This group is $R_n$, the group of permutations on the pairs of vertices induced by $S_n$, where $S_n$ is the usual group of permutations on $n$ symbols ($n$ vertices of a cgraph).

**Example 2.6.1:** Let $n$=3, the field be $GF(3) = \{0,1,2\}$. The colors be white (= 0), red (= 1) and blue (= 2). The figure counting series for this case will be $A(x, y) = 1 + x + y$. The cycle index

$$Z(R_3) = \frac{1}{6}(t_1^3 + 3t_1t_2 + 2t_3)$$

Putting $1 + x + y$ for $t_1$, $1 + x^2 + y^2$ for $t_2$ and $1 + x^3 + y^3$ for $t_3$, we get the following configuration counting series:

$$B(x) = 1 + x + y + x^2 + y^2 + xy + x^2y + xy^2 + x^3 + y^3.$$

Thus, there are in all 10 types of colored graphs, colored in say red and blue color and white color indicating absence of edge, on three points and they are as follows:

(1) A cgraph containing three white lines (totally disconnected graph containing three vertices in usual graph theoretic sense).
(2) A cgraph containing one red line and two white lines.
(3) A cgraph containing one blue line and two white lines.
(4) A cgraph containing two red lines and one white line.
(5) A cgraph containing one red line, one blue line, and one white line.
(6) A cgraph containing two blue lines and one white line.
(7) A cgraph containing one red line and two blue lines
(8) A cgraph containing one blue line and two red lines
(9) A cgraph containing three red lines.
(10) A cgraph containing three blue lines.

**3. New Preliminaries:** In order to get some interesting results special for the cgraphs we should work in Galois filed modulo ($n$) for some prime $n$ and exploit the special constructs that follow from it.

Some definitions are in order:

**Definition 3.1:** A **cgraph** *G*, containing *m* vertices and *n* colored edges consists of a vertex set $V(G) = \{v_1, v_2, \cdots, v_m\}$ and an edge set $E(G) = \{e_1^{C_{j1}}, e_2^{C_{j2}}, \cdots, e_n^{C_{jn}}\}$, where the superscript $C_{j_r}$ on edge label $e_r$ (as expressed usually) indicates that the edge $e_r$ is colored with color $C_{j_r}$.

When there is no confusion, one can drop the symbol *C* in the edge labeling and simply write the edge set as $E(G) = \{e_1^{j1}, e_2^{j2}, \cdots, e_n^{jn}\}$, where numbers $\{j_1, j_{2,} \cdots, j_n\}$ represent the corresponding colors, and they are essentially elements of the Galois field modulo (*n*).

**Definition 3.2:** A **cdigraph** is a directed cgraph.

**Definition 2.3:** The $\pi$ **-Complement** of a cgraph $G, \pi(G)$, is a cgraph obtained by application of a permutation on *n* symbols {0, 1, 2, ⋯, *n*—1} representing the color change of the edges of the given cgraph $G$. This leads to change in the color of every edge (or, no color change for the edge with color *j* if $j \to j$ under the action of $\pi$) of the cgraph $G$ having color *j* to a new color $\pi(j)$ in the complement cgraph $\pi(G)$.

**Definition 3.4:** Two cgraphs $G$ and $H$ are **complements** of each other if $H$ is $\pi$ -Complement of $G$ for some permutation $\pi \in S_n$.

**Proposition 3.1:** $S_n$, the group of permutations on *n* symbols {0, 1, 2, ⋯, *n*—1} forms the **complementomorphism group** for any cgraph.

**Proposition 3.1:** $G \cong H$ if and only if $\pi(G) \cong \pi(H)$.

**Proof:** By isomorphism of *G* and *H*, if there is a *k*-colored edge between $v_i, v_j \in G$ then there will be a *k*-colored edge between $\phi(v_i), \phi(v_j) \in H$, where $\phi$ is the corresponding isomorphism map. Now, the existence of *k*-colored edge between $v_i, v_j \in G$ implies the existence of $\pi(k)$-colored edge between $\pi(v_i), \pi(v_j) \in \pi(G)$ which in tern implies the existence of $\pi(k)$-colored edge between $\pi(\phi(v_i)), \pi(\phi(v_j)) \in \pi(H)$. Converse is similar and follows from the fact that $\pi^{-1}$ exists.

**Definition 3.5:** A **subcgraph** $H$ of a cgraph $G$ is a cgraph such that $V(H) \subseteq V(G)$ and $E(H) \subseteq E(G)$. The **subcgraph** $H$ is called **induced** if it is a maximal csubgraph.

**Definition 3.6:** A **cgraph** is called ***j*-complete** if there exists an edge joining every vertex pair of $G$ and all these edges have color *j*.

**Definition 3.7:** A subset $S \subseteq V(G)$ is called ***k*-independent set** if the induced subcgraph induced by $S$, say $G[S]$, has no edge having color *k*.

Note that a vertex set of a *j*-complete cgraph forms a *k*-independent set for all $k \neq j$.

**Definition 3.8:** A cgraph is called ***k*-bipartite** if its vertex set is union of two disjoint sets such that there are edges (at least one) connecting vertices in the disjoint partitioning having color *k* and having one end point from each disjoint set.

**Definition 3.9:** A ***k*-complete bipartite** cgraph is a *k*-bipartite cgraph having an edge set made up of all vertex-pairs having a vertex from each of two disjoint sets and all these edges are colored by *k* color.

**Definition 3.10:** A ***k*-path** in a cgraph is a cgraph as a single vertex or an ordered list of distinct vertices $v_1, v_2, \cdots, v_n$ such that $(v_j, v_{(j+1)})$ is a *k*-colored edge for all $1 \leq j \leq (n-1)$.

**Definition 3.11:** A ***k*-cycle** is a *k*-path which is closed, i.e. $v_n v_1$ is also a *k*-colored edge.

**Definition 3.12:** A cgraph $G$ is called ***j*-connected** if it has a *j*-colored path joining each vertex pair $(v_i, v_k)$ such that $v_i, v_k \in V(G)$.

Note that a *j*-connected cgraph is the one which contains a connected spanning subcgraph (subctree) having all edges colored with color *j*.

**4. Matrices $A(G)$ and $I(G)$ and cisomorphism of cgraphs:** In order to specify a cgraph we need idea of cisomorphism.

**Definition 4.1:** Two cgraphs $G$ and $H$ are **cisomorphic** if there exists a bijection (called cisomorphism) $\phi : V(G) \to V(H)$ such that if there is a *j*-colored edge joining vertices $v_i, v_k \in V(G)$ then there exists a *j*-colored edge joining vertices $\phi(v_i), \phi(v_k) \in V(H)$.

Note that cisomorphism preserves both, the adjacency and color.

**Definition 4.2:** Given a cgraph $G$ with vertex set $V(G) = \{v_1, v_2, \cdots, v_m\}$ and edge set $E(G) = \{e_1^{j1}, e_2^{j2}, \cdots, e_n^{jn}\}$, then the **adjacency matrix** of $G$ is a matrix $A(G) = [a_{ij}]$, of size $m \times m$ such that $a_{ii} = 0$, indicating absence of self-loops, and $a_{ij} = k$, when $i \neq j$, if there exists a ***k*-colored** edge joining vertices $v_i, v_j$ of $G$, $k \in \{0, 1, 2, \cdots, n-1\}$.

**Theorem 4.1:** A cgraph $G$ is made up of *j*-colored subcgraphs $G_j$ and adjacency matrix of cgraph $G$, $A(G)$, can be expressed as sum of adjacency matrices of *j*-colored subcgraphs, $A(G_j)$, where $A(G_j)$ is obtained by replacing those entries of $A(G)$ which are not equal to *j* by 0 (Note that $j + 0 = 0 + j = j$).

**Definition 4.3:** Given a cdigraph $G$ with vertex set $V(G) = \{v_1, v_2, \cdots, v_m\}$, and edge set of directed edges $E(G) = \{e_1^{j1}, e_2^{j2}, \cdots, e_n^{jn}\}$, then the **adjacency matrix** of $G$ is a matrix $A(G) = [a_{ij}]$, of size $m \times m$ such that $a_{ij} = k$, if there exists a ***k*-colored directed** edge joining vertices $v_i$ to $v_j$ of $G$, $k \in \{0, 1, 2, \cdots, n-1\}$, and $a_{ij} = 0$, if edge $(v_i, v_j)$ is absent (or white colored one).

**Proposition 4.1:** Two cgraphs $G$ and $H$ are cisomorphic if and only if there exists a permutation $\sigma$ with permutation matrix $\mathbf{M}$ such that $A(G) = \mathbf{M}\, A(H)\, \mathbf{M}^{\mathbf{T}}$, where $\mathbf{M}^{\mathbf{T}}$ is the transpose of $\mathbf{M}$.

**Proposition 4.2:** Cisomorphism is an equivalence relation.

**Definition 4.4:** A **cisomorphism class** of cgraphs is an equivalence class of cgraphs under the cisomorphism relation.

**Definition 4.5:** By an **unlabeled** (i.e. no labels to vertices) cgraph we meant a cisomorphism class of cgraphs.

**Definition 4.6:** A **cautomorphism** of a cgraph $G$ is a permutation of $V(G)$ that is cisomorphism from $G$ to $G$.

**5. Reconstruction Conjectures for Cgraphs:** The vertex reconstruction conjecture for cgraphs can be stated as follows:

**Definition 5.1:** Two cgraphs with at least three vertices $G$ and $H$ are called **chypomorphic** if there exists a bijection called **chypomorphism**
$\sigma : V(G) \to V(H)$
$v \to u$, such that
$G - v$ is cisomorphic to $H - u$.

**Conjecture 5.1 (Reconstruction Conjecture):** Chypomorphic cgraphs are cisomorphic.

**Conjecture 5.2 (Edge Reconstruction Conjecture):** If $G$ and $H$ are two cgraphs with at least four edges, and if there exists a bijection $\theta : E(G) \to E(H)$ such that $G$−e is cisomorphic to $H - \theta$ (e) then $G$ and $H$ are cisomorphic.

**6. Vector Spaces and Subspaces Associated with cgraphs:** It is well-known fact that $F = GF(p)$, forms a finite field containing elements {0, 2, 3, …, $p-1$} if and only if $p$ is a prime number. This field is called Galois field modulo $p$. We can associate a vector space basis made up of vectors: {$x_1$= (1, 0, 0, …, 0), $x_2$ = (0, 1, 0, …, 0),…, $x_q$ = (0, 0, 0, …, 1)} with cgraphs. Every ($p$, $q$) cgraph can be expressed as a vector

$x_G = \sum_{j=1}^{q} c_j x_j$, where $c_j \in GF(p)$, the Galois field modulo ($p$), $p$ prime, and the number $c_j$ stands for the color $C_j$ .Thus, the cgraph can be represented by vector $(c_1, c_2, \cdots, c_q)$ where all $c_j \in GF(p)$.

We can associate, in all, $p^q$ vectors, representing $p^q$ cgraphs with a given cgraph represented by vector $(c_1, c_2, \cdots, c_q)$ as follows: (1) There will be cgraphs represented by vectors $(j_1, j_2, \cdots, j_q)$ such that $0 \le j_m \le c_m$, such that $0 \le m \le q$ called **subcgraphs**.
(2) There will be cgraphs represented by vectors $(k_1, k_2, \cdots, k_q)$ such that $c_m \le k_m \le p-1$, such that $0 \le m \le q$ called **supercgraphs**.
(3) There will be other cgraphs which partially satisfy condition in (1) and partially satisfy condition in (2) and so they are **neither subcgraphs nor supercgraphs** of the given cgraph.

Thus, there exists a vector space, $W_G$, associated with each cgraph $G$, and this vector space consists of
(1) Finite field, $F = GF(p)$.
(2) $p^q$ number of $q$-tuples where $q$ is the total number of colored edges in $G$.
(3) An addition operation between two vectors $\mu, \nu$ in this space, defined as vector sum, denoted by operation $\oplus$. Thus, let $\mu = (\mu_1, \mu_2, \cdots, \mu_q)$ and $\nu = (\nu_1, \nu_2, \cdots, \nu_q)$ then $\mu \oplus \nu = (\mu_1 + \nu_1, \mu_2 + \nu_2, \cdots, \mu_q + \nu_q)$, where "+" being the addition modulo($p$).
(4) Scalar multiplication between scalar $c \in GF(p)$ and a vector $\mu$ defined as $c \bullet \mu = (c\mu_1, c\mu_2, \cdots, c\mu_q)$.
Note that here the zero vector, the $q$-tuple of zeros, forms the additive identity, i.e. the identity for vector sum operation.

We now conclude with the following hope:

**Conclusion:** It is likely that the readers will find many more interesting results and thus carry out the program proposed in this small paper and also achieve the essential development of the ideas required for the colorful graph theory in much more meaningful way.

## References

1. Mehendale Dhananjay P., Finite Projective Planes, http://arXiv.org/abs/math/0611492